# Evaluating Stochastic Methods in Power System Operations with Wind Power

Yishen Wang
University of Washington
Seattle, WA, USA
Email: ywang11@uw.edu

Zhi Zhou, Cong Liu, Audun Botterud
Argonne National Laboratory
Argonne, IL, USA
Email: {zzhou, liuc, abotterud}@anl.gov

*Abstract*—Wind power is playing an increasingly important role in electricity markets. However, it's inherent variability and uncertainty cause operational challenges and costs as more operating reserves are needed to maintain system reliability. Several operational strategies have been proposed to address these challenges, including advanced probabilistic wind forecasting techniques, dynamic operating reserves, and various unit commitment (UC) and economic dispatch (ED) strategies under uncertainty. This paper presents a consistent framework to evaluate different operational strategies in power system operations with renewable energy. We use conditional Kernel Density Estimation (KDE) for probabilistic wind power forecasting. Forecast scenarios are generated considering spatio-temporal correlations, and further reduced to lower the computational burden. Scenario-based stochastic programming with different decomposition techniques and interval optimization are tested to examine economic, reliability, and computational performance compared to deterministic UC/ED benchmarks. We present numerical results for a modified IEEE-118 bus system with realistic system load and wind data.

## I. INTRODUCTION

The rapidly growing penetration of renewable energy resources has created significant challenges for the electricity markets and grid operators due to its inherent variability and uncertainty. Commonly adopted approaches to address these challenges include advanced forecasting techniques [1-3] and improved system operational uncertainty modeling [4-16]. Traditionally, the system operators make the commitment and dispatch decisions based on a deterministic renewable point forecast. Numerous statistical methods and machine learning techniques have been applied to the accuracy of wind power forecasting [1-3], including Gaussian processes, support vector machines, artificial neural networks, random forest (RF), etc. [1]. Although these methods produce a point forecast with improved accuracy, additional information about the uncertainty in the forecast is needed. Recently, probabilistic wind power forecasting has attracted substantial attention as a means to provide wind power's probability distribution in addition to the central point forecast [2]. Kernel density estimation [3], quantile regression and other parametric or non-parametric algorithms [1-2] have been proposed in this context. With such probabilistic information, wind power scenarios or uncertainty sets with a certain confidence level can be properly constructed to feed power system operational models for improved representation of wind power forecast uncertainty.

In the literature, many researchers have applied different operational schemes to incorporate stochastic wind power generation [4-16]. These approaches mainly includes deterministic unit commitment (DUC) with dynamic reserves [4-6], scenario-based stochastic UC (SUC) [4-7], interval-based UC (IUC) [8-10], robust UC (RUC) [11-12], chance-constrained UC [13], as well as hybrid and unified approaches [14] [15]. The traditional DUC typically keep contingency reserves to address the possibility for a generation or line outage. To better account for the stochastic system conditions, including forecast uncertainty, more dynamic operating reserve strategies can be used. For instance, Zhou et al [4] propose that the dynamic reserve (i.e. reserve for variable generation) amount is determined through the difference between a central point forecast and a lower level quantile obtained from a probabilistic forecast. In contrast, scenario-based SUC determines the unit commitment decisions to provide implicit reserves by considering a set of forecast scenarios. Papavasiliou et al [6] compare the two-stage SUC and DUC with different reserve policies. The results show that SUC can produce a lower cost solution than DUC. This conclusion is similar to the one drawn in [4]. In [7], Uckun et al. propose an improved stochastic UC formulation which approximates a multi-stage formulation. The results show a significantly lower operational cost than traditional two-stage SUC, with limited increase in computational burden. SUC requires a good quality scenario set which fully characterizes the wind uncertainty. Sub-hourly SUC is discussed in [29]. However, SUC is generally not computationally tractable with more than 100 scenarios for a system of realistic size. Hence, scenario reduction is applied to select representative scenarios [16]. Fast forward selection performs well for both economical and computational aspects [17]. Decomposition frameworks are also helpful to address the computational challenges. We consider Benders decomposition (BD) and linear-shift-factor (LSF) decomposition [23], two common approaches to improve the computational performances in UC, as discussed in section III.

IUC does not require full information about the wind power distributions as it only uses the upper and lower bound for the wind power forecast. Ramping capabilities between these bounds are considered to ensure the operational feasibility in extreme situations [8]. Pandzic et al [9] modify the IUC approach by determining the ramping requirements with a data-driven method to reduce the over-conservativeness within the original formulation. Liu et al [10] further improve the IUC formulation with a fuzzy set approach to transform the crisp interval bounds into a fuzzy one, and shows that if wind curtailment cost is zero, then only the lower bound is necessary. Similarly, RUC only requires a polyhedral uncertainty set. It aims to minimize the worst case operating cost and ensures the feasibility of all possible realizations within the set [11-12]. Hence, the RUC approach is by definition a very conservative operational strategy.

The contribution of this paper is to systematically evaluate the impact of UC formulations, operating reserve policies, and wind uncertainty representations on short term power system operations in terms of economic, reliability, and computational performance. We focus the comparison on DUC, IUC, and SUC formulations with different decomposition techniques. The paper is organized as follows. Section II describes the probabilistic forecasting techniques and scenario generation and reduction methods used in the analysis. Section III briefly introduces the different UC formulations and decomposition techniques that are compared in the case study. The market simulation platform is presented in section IV, and section V presents the case study based on a modified IEEE-118 system. Conclusions are summarized and discussed in section VI.

## II. PROBABILISTIC WIND POWER FORECASTING AND SCENARIO GENERATION

We adopt a time-adaptive quantile-copula to produce the probabilistic wind power forecast. Following the approach in [3], the point forecast and the hour of the day are used as explanatory variables in the model, which is trained on historical wind power data. Then, conditional Kernel Density Estimation (KDE) is applied to estimate the conditional wind power probability distribution. Next, the time-adaptive quantile-copula approach estimates the forecast probability density function for each time interval. We refer to [3] for more details. From the probabilistic forecasts we generate the wind uncertainty sets based on desired confidence levels for the DUC and IUC formulations. For instance, for a 90% confidence range, we use the 5% and 95% quantiles as the lower and upper bounds for wind power. To generate scenarios for multiple wind farm for the SUC formulation, a spatio-temporal covariance matrix is constructed following [18]. Then, a large set of scenarios is sampled from the estimated probability distribution and covariance matrix. Finally, the scenario reduction [16] is applied to obtain a reduced scenario set with GAMS/SCENRED [19]. An example of a reduced scenario set is shown in Figure 1.

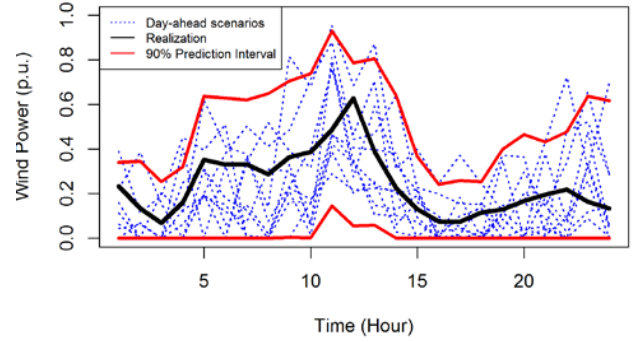

Figure 1 Reduced scenarios and confidence interval for wind power

## III. UC FORMULATION AND DECOMPOSITION

Below we present high-level formulations of DUC, IUC, and SUC with extensive form (SUC-E), SUC with linear-shift-factor decomposition (SUC-LSF) and SUC with Benders decomposition (SUC-BD). The formulations are presented in abstract matrix format.

### A. Deterministic Unit Commitment with Dynamic Reserve

DUC is the current industry practice and there are many papers devoted to improve the computational issues for DUC formulations. We adopt the tight and compact formulation proposed by Morales-Espana in [20]. We model multiple reserve products including the regulation up/down, spinning and non-spinning reserve. The compact DUC formulation is:

$$\min \quad C^{su}u + C^{fuel}p \quad (1)$$
$$Au \leq B \quad (2)$$
$$Cu + Dp \leq M \quad (3)$$
$$Fu + Gp \leq N \quad (4)$$

where $u$ represents the commitment, start-up and shut-down binary variables and $p$ denotes the dispatch related continuous variables. In (1) the model minimizes the sum of startup cost and fuel costs. Equation (2) represents the constraints related to the binary variables, including minimum up and down times and commitment status. Equation (3) represents all other constraints of the UC model, including minimum/maximum generation capacity, wind maximum limit, ramping constraints, power balance, and line flow limits [20]. Equation (4) represents all the reserve constraints. For more detailed formulations we refer to [4-5] and [20].

The dynamic reserve amount is determined from the probabilistic wind power forecast, normally we choose the 10%, 5% or 1% quantile to hedge against over-forecasting. The sum of this dynamic reserve and a fixed contingency reserve satisfy the system wide total operating reserve requirements. Regulation up/down, spinning and non-spinning reserve requirements are determined as fixed fractions of the total reserve requirements. The percentages are obtained from PJM historical data statistics [21].

## B. Interval Unit Commitment

The IUC formulation proposed in [10] is adopted in this paper. In addition to equations (1)-(3), equation (5) is added to ensure system feasibility when the wind lower bound is hit. This includes the power balance, ramping constraints, reserve limits and line flow limits in the lower bound and in all transitions to the lower bound. The base contingency reserve in equation (4) is still included in the formulation. Since the IUC optimization determines the necessary additional reserve amount to account for wind uncertainty, the explicit dynamic reserve is no longer needed in this case. The objective function is to minimize the total costs for the central point forecast, same as in the original DUC case.

$$Ju + Kp \leq L \quad (5)$$

The detailed IUC model is presented in [10] with proof of the formulation and its fuzzy counterparts. Noted we are using a crisp lower bound here, and not a fuzzy constraint as presented in [10].

## C. Stochastic Unit Commitment

We consider a two-stage SUC formulation in this paper:

$$\min \quad C^{su}u + E[C^{fuel}p(\xi)] \quad (6)$$
$$Au \leq B \quad (7)$$
$$Cu + Dp(\xi) \leq M(\xi) \quad (8)$$
$$Fu + Gp(\xi) \leq N(\xi) \quad (9)$$

The two-stage SUC minimizes the expected operational cost with respect to the dispatch related constraints for each scenario denoted with $\xi$. Note that equation (9) is added here to enable the use of a small amount of additional reserve for each scenario to address uncertainty not included in the scenario set and further improve the system performances as presented in [4].

## D. LSF Decomposition

When the system is large, it is computationally expensive to evaluate the power flow limits for all lines. As shown in [23], the number of congested lines during operations is typically only a small fraction of the full line set. Therefore, only considering critical line flow limits can improve computational performances without sacrificing solution accuracy. This assumption also reflects real-world ISO practice: In [22], only a selected number of lines are labeled with "Monitor Flag" to check their line limit during operation, and other line limits are not evaluated in the constraints.

In the SUC problem, considering the number of lines, hours and scenarios, only a few transmission constraints are actually binding. Therefore, applying LSF decomposition can effectively relieve the computational burden as checking flow limits is essentially a LP feasibility problem, which can be solved very efficiently. The pseudo code is as follows:

**Algorithm 1 - LSF Decomposition:**
1) *Initialize the line checking set S (l,t,ξ) as an empty set, l, t, ξ is the index for line, time, scenario respectively;*
2) *Solve the SUC problem with network constrints only imposed on the line checking set S;*
3) *Evaluate all the remaining line tuples (l, t, ξ) for their flow limits: if no limits is violated, go to step 4); else, add the violated line tuples (l, t, ξ) in S, and go to step 2);*
4) *Stop the program, and output the result.*

## E. Benders Decomposition

For the two-stage SUC structure, as the second stage economic dispatch problem per scenario is an LP problem, the BD can be applied here. The SUC problem is suitable to decompose into a master problem and a collection of subproblems for each scenario. The master problem involves all the UC constraints and Benders cuts. In the classical Benders approach, the subproblems, with the dispatch constraints [23], are used to generate feasibility and optimality cuts for the master problem to converge. However, the convergence rate can be slow. Similar to the approach in [24], we only consider the optimality subproblems, and slack variables are added to form complete second-stage recourse decisions. In each iteration, the subproblems evaluate the dispatch decision based on commitment status sent from the master problem, and return sensitivities to generate Benders cuts for the master solution. In the master problem, we add an online generation capacity constraint and an aggregate optimality cut summing all optimality cuts together to further improve the convergence. The pseudo code is provided below, with more details provided in [23-24].

**Algorithm 2 - Benders Decomposition:**
1) *Initialize the problem, UB=+∞, LB=-∞;*
2) *Solve master problem, update the LB with objective function $LB=z_{LB}$;*
3) *Run optimality subproblem for each scenario: update the upper bound with solution from optimality subproblems $UB = z_{UB}$;*
4) *Check the convergence: if |UB-LB|≤ε, go to step 5); else, generate optimality cuts and return to step 2);*
5) *Stop the program, and output the results.*

Note that the above algorithm follows a modified BD approach, and multiple acceleration techniques can be used to further improve the speed of convergence [23-24].

## IV. ELECTRICITY MARKET FRAMEWORK

The simulation is based on a rolling decision-making process typically applied in U.S. electricity markets. This includes the Day-ahead UC (DA-UC), Reliability Assessment Commitment (RAC), and Real-Time Economic Dispatch (RT-ED), with updated forecast information applied at each stage. In DA-UC, the unit commitment schedule for the day-ahead is determined and passed to the RAC stage. In RAC, the slow generators follow the DA schedule, whereas quick-start generators are allowed to adjust their status based on the updated forecast. This schedule update passes to the RT-ED for the actual dispatch solution. This final dispatch then passes to the next day's simulation as the initial condition. DA-UC and RAC considers full reserve constraints based on the forecast and scenarios. In contrast, RT-ED only procures the contingency reserve. The flowchart for the simulation is presented in Figure 2. A four month out-of-sample simulation is carried out for the different operational strategies to cover daily, weekly, and seasonal patterns in loads and wind resource availability.

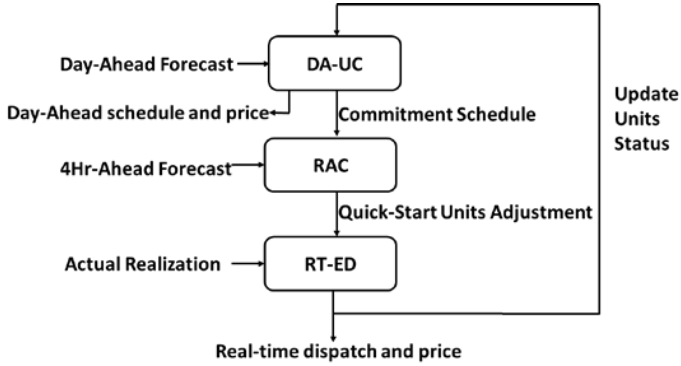

Figure 2 Rolling Electricity Market Framework

## V. NUMERICAL RESULTS

### A. Test System

A modified IEEE 118-bus system [25] is applied as the test system. There are a total of 186 lines and 54 generators, including 21 quick-start units. The load is scaled based on actual load profiles from the U.S. state of Illinois in 2006. The system peak load is 6616 MW and the average load is 4093 MW in the four month simulation period from July to October. The energy and reserve cost data are derived from the original IEEE test system and PJM data [21, 25]. There are three wind farms located at bus 15, 54, 96, i.e. one in each zone. The installed wind capacity is 902.6, 911.4, and 1232.2 MW respectively. The average available wind energy is 21% of the load. The wind profiles are selected from the Eastern Wind Integration Dataset [26]. The Day-ahead forecast is used as the initial forecast for the DA-UC stage, the four-hour-ahead forecast is used as the forecast for the RAC stage. Out-of-sample, actual wind realizations are applied at the RT-ED stage to evaluate the UC performances. Probabilistic forecast for the three wind farms are first generated, then 1000 scenarios are sampled each day considering spatio-temporal correlations, with a reduced scenario set of 10 selected, as explained in section II. We use 10 scenarios primarily to limit the computational burden, especially for the extensive form. However, the extra economic benefits of more scenarios are also limited, according to [17]. In the case study, regulation up/down equals to 15% of the total reserve requirement; while spinning and non-spinning reserve each takes 42.5% of the total requirement.

The framework is implemented in AMPL [27] and the MIP solver is IBM ILOG CPLEX 12.1 with 0.1% optimality gap [28]. All the results are obtained on a Linux-based server at Argonne National Laboratory with 2x Intel E5430 Xeon CPUs and 32 GB of RAM.

### B. Case Study Description

The DUC with perfect wind forecast (DUC-perfect) and point forecast (DUC-point) are selected as the two reference cases. The point forecast is obtained as the 50% quantile from the probabilistic wind power forecast. In these two reference cases, only contingency reserve is imposed at all stages. For DUC with dynamic reserve (DUC-DR), we choose the difference between 50% and 10% quantiles as the dynamic reserve quantity. Correspondingly, the 10% quantile is also chosen as the lower bound for the IUC formulation. Sensitivity analysis showed that using the 10% quantile is more economical than using lower quantiles, as the benefits of the dynamic reserve cannot compensate for the additional reserve costs. In the SUC cases, in addition to the different solution techniques, we compared different levels of additional wind reserve (0, 1%, 5% and 10% of wind power at each scenario). In Tables 1 and 2, we only present the result with the best performance, i.e. with 0 wind reserve.

### C. Economic and Reliability Performances

We present economic and reliability results for summer (July and August) and fall (September and October) in Table 1. DUC-perfect always produces the cheapest solution at all times, as expected given the perfect forecast assumption.

Table 1 Average Daily Cost and Load Shedding for summer and fall

| UC Model | Summer | | | Fall | | |
|---|---|---|---|---|---|---|
| | DA Cost ($M) | RT Cost ($M) | RT Load Shed (MW) | DA Cost ($M) | RT Cost ($M) | RT Load Shed (MW) |
| DUC-perfect | 1.517 | 1.517 | 0 | 1.023 | 1.023 | 0 |
| DUC-point | 1.589 | 1.523 | 0 | 1.074 | 1.050 | 3.43 |
| DUC-DR-10 | 1.637 | 1.523 | 0 | 1.146 | 1.043 | 0.52 |
| IUC-10 | 1.623 | 1.524 | 0 | 1.131 | 1.046 | 0 |
| SUC-E-0 | 1.565 | 1.521 | 0 | 1.078 | 1.034 | 0 |

In the summer, net load is relatively large, more generators are on, and online capacity is sufficient to mitigate the forecast inaccuracy from DA scheduling to RT-ED. In this case, when the system buys more explicit reserve as in DUC-DR and schedules more implicit reserve in IUC, the additional cost from this reserve component does not bring economic benefits. In contrast, when the net load is relatively lower in the fall, online capacity is not enough to handle the extreme events, and load shedding happens at particular days under DUC-point. In these cases, dynamic reserve and IUC produce a more reliable schedule, and the additional reserves result in a cost saving due to lower load curtailment. In sensitivity runs, when we choose a more conservative strategy and select a lower quantile, as 5% and 1%, the reliability for IUC holds, and for dynamic reserve the load shedding is reduced. However, as we compare the cost components, dynamic reserve pays more for buying reserve in the DA from slow generators and produces a cheaper schedule in RT. Whereas, IUC provide more implicit reserve from quick-start units in DA with less reserve bought, but their schedule is more expensive in RT. For SUC cases, the scenarios capture the wind uncertainty better. In the DA-UC stage, SUC encourages more flexible units (including quick-start units and slow-start units with large ramping capabilities) online and the wind utilization rate is higher with less wind spillage, which results in a similar range of DA costs compared to DUC-point. As the wind uncertainty is better managed with SUC, all the SUC RT-ED solutions are very close to the DUC-perfect and cheaper than other DUC and IUC cases, both for summer and fall. At the same time, there is no load shedding at all showing that the optimized schedule is also reliable. Hence, the results indicate that SUC reaches a better balance between system costs and reliability.

In Table 2, we present the average daily number of commitments for fast and slow generators. As the two reference DUC strategies (perfect and point) only rely on a point forecast, they tend to use the cheaper slow units as much as possible. Instead, for the dynamic reserve and IUC cases, they consider the forecast uncertainty and more quick-start units are running for these cases. As IUC provides implicit reserve and ensure the reserve deliverability for the transition and lower bound cases, they tend to start up more quick-start units with a more expensive schedule. The day-ahead cost saving from IUC compared to dynamic reserve is mainly from the avoidance of reserve curtailments penalties. The results also show that the stochastic cases commit fewer units to hedge against uncertainty. The quick-start units only run online when it is necessary. Therefore, their operational schedules are close to the two reference cases, but with better reliability metrics than other operational strategies (Table 1). Table 1 and Table 2 both illustrate the advantage of stochastic scheduling methods.

Table 2 Average Daily Total Commitments

| Average Daily Number of Commitments | Summer | | | Fall | | |
|---|---|---|---|---|---|---|
| | Fast Units | Slow Units | Total | Fast Units | Slow Units | Total |
| DUC-perfect | 1 | 706 | 707 | 2 | 564 | 566 |
| DUC-point | 2 | 728 | 730 | 1 | 597 | 598 |
| DUC-DR -10 | 5 | 733 | 738 | 20 | 627 | 647 |
| IUC-10 | 23 | 761 | 784 | 30 | 707 | 737 |
| SUC-E-0 | 4 | 711 | 715 | 1 | 567 | 568 |

In sensitivity analysis, as we compare the SUC cases with different additional scenario reserves, this dynamic reserve component sometimes brings economic benefits in RT. The reason is that the additional reserve slightly increases the total number of online flexible units, which can react better in real-time. This conclusion is also similar to the findings in [4]. Moreover, different decomposition generally gives very similar results, with differences between solutions all within the optimality gap. The advantage from decomposition is more from the computational aspect, as we will discuss next.

*D. Computational Performances*

As shown in Table 3, compared to the DUC with point forecast, procuring dynamic reserve increase the computing time, but the overall speed is still very fast. Interestingly, IUC outperforms DUC from a computational perspective, and have a very similar performance with the two DUC reference cases. The reason is that IUC add a set of feasibility check constraints to ensure the system has enough reserve, and these constraints drive the system to schedule implicit reserve only when it is necessary, which results in less dynamic reserve scheduled in the IUC case. Therefore, IUC adds less computational burden than the dynamic reserve case.

Table 3 Average Computational Time Per Day

| UC Model | DA Time(s) | RAC Time (s) | RT Time (s) |
|---|---|---|---|
| DUC-perfect | 11 | 4.3 | 0.54 |
| DUC-point | 10 | 4.6 | 0.89 |
| DUC-DR-10 | 22 | 5.9 | 1.01 |
| IUC-10 | 11 | 4.3 | 0.46 |
| SUC-E-0 | 248* | 62 | 1.19 |
| SUC-LSF-0 | 266 | 48 | 0.87 |
| SUC-BD-0 | 171 | 14 | 0.89 |

\* Out of memory issues for some days.

As expected, the computation time is significantly longer for the SUC cases compared to IUC and DUC. Still, with ten scenarios, the times are all below 10 minutes for the day-ahead problem. Note that although the extensive form has similar computing times as the LSF decomposition, there are out-of-memory issues for particular days where AMPL fails to provide a solution. Hence, the optimality gap and other settings have to be relaxed and adjusted to provide a solution. BD is faster than extensive form and LSF decomposition, but some efforts are needed to construct proper subproblems for better convergence. In contrast, with LSF decomposition, these extra manual efforts are not necessary. In general, decomposition will allow for incorporating more scenarios, something we plan to test in future work.

Table 4 shows that LSF decomposition generally converges in 3 to 4 iterations. Each cut is essentially representing one congested or potentially congested line. On average only 0.73% of the lines are potentially congested for this case. A better initialization step could place pre-defined critical lines in the line checking set to further reduce the computing time. Finally, Table 4 also shows that BD can normally converge in

15 iterations to reach the optimality criterion for the day-ahead problem. However, for some extreme cases the maximum iteration limit of 100 is met. In contrast, after fixing the slow generators, RAC only needs one more iteration to find the optimal solution.

Table 4 Decomposition Convergence Statistics

| UC Model | DA Average /Maximum Number of Iteration | DA Average /Maximum Number of Cuts | RAC Average /Maximum Number of Iteration | RAC Average /Maximum Number of Cuts |
|---|---|---|---|---|
| SUC-LSF-0 | 3.7/7 | 328.0/940 | 3.4/4 | 293.2/702 |
| SUC-BD-0 | 15/100 | 150/1100 | 2/2 | 11/11 |

## VI. CONCLUSION

As more wind power is being integrated into the electric power grid, the forecast uncertainty brings operational challenges for the system operators. In this paper, different operational strategies for uncertainty management are presented and evaluated. A consistent simulation framework is used to analyze the performance of different reserve policies and scheduling techniques. The numerical results from the IEEE 118-bus system show that SUC formulations provide a reliable schedule without high increases in costs. Moreover, the results indicate that decomposition can help in overcoming the computational obstacles for SUC. In contrast, DUC and IUC tend to give higher system costs as more reserves are being scheduled, but require much lower computational time.

## VII. ACKNOWLEDGMENT

This work was supported by the U.S. Department of Energy through its Wind and Water Power Program. The submitted manuscript has been created by UChicago Argonne, LLC, Operator of Argonne National Laboratory ("Argonne"). Argonne, a U.S. Department of Energy Office of Science laboratory, is operated under Contract No. DE-AC02–06CH11357.